\documentclass[letterpaper, 10 pt, conference]{ieeeconf}  

\IEEEoverridecommandlockouts                              

\usepackage{comment}
\usepackage{cite}
\usepackage{amsmath,amssymb,amsfonts}

\DeclareMathOperator{\sgn}{sgn}
\usepackage{algorithmic}
\usepackage{graphicx}
\usepackage{textcomp}
\usepackage{xcolor}
\usepackage[normalem]{ulem}
\def\BibTeX{{\rm B\kern-.05em{\sc i\kern-.025em b}\kern-.08em
    T\kern-.1667em\lower.7ex\hbox{E}\kern-.125emX}}

\newtheorem{theorem}{Theorem}
\newtheorem{corollary}{Corollary}

\newtheorem{remark}{Remark}

\usepackage{microtype}

\title{\LARGE \bf \spaceskip=1\fontdimen2[-30]{From Hamilton--Jacobi--Bellman to Riccati via Carleman: The Scalar Case}
}

\author{Philip J. Elias$^{1}$, Qiyu Sun$^{2}$, and  Nader Motee$^{3}$ 
\thanks{*This work was supported in part by the ONR under grants numbers
N00014-23-1-2779 and N00014-26-1-2246.}
\thanks{$^{1}$P. J. Elias is with the Department of Chemical and Biomolecular Engineering, Lehigh University, Bethlehem, PA 18015, USA
        {\tt\small phe222@lehigh.edu}}%
\thanks{$^{2}$Q. Sun is with the Dr. Yvette Kanouff School of Data, Mathematical and Statistical Sciences,  University of Central Florida, Orlando, FL 32816, USA. 
        {\tt\small qiyu.sun@ucf.edu}}%
\thanks{$^{3}$N. Motee is with the Department of Mechanical Engineering and Mechanics, Lehigh University, Bethlehem, PA 18015, USA
        {\tt\small motee@lehigh.edu}}%
}

\begin{document}

\maketitle
\thispagestyle{empty}
\pagestyle{empty}

\begin{abstract}
The Hamilton--Jacobi--Bellman (HJB) equation of a nonlinear system is a
nonlinear partial differential equation, and in that form it gives little
away. This paper shows that it hides a Riccati equation. Carleman
linearization reveals this structure. For scalar analytic input-affine systems, the lifted HJB equation is a single infinite-dimensional
Riccati equation in a lower-triangular Toeplitz matrix. Its first entry is
the Riccati equation of the standard linearization at equilibrium. Its
remaining entries hold the nonlinear terms. Nothing about this is an approximation.
 Under controllability of the standard linearization and detectability of the
cost, the equation has exactly one solution with positive diagonal entries,
and it can be computed one entry at a time. This solution is the derivative of the value
function. We prove that finite-sections of the solution are exact and converge exponentially to
the derivative of the value function near the origin. The scalar case has a closed-form solution. This closed form is our benchmark: every result on the lifted equation is checked against it. Nothing in this structure depends on the state dimension. Our methodology points to a general theory of optimal control in which the HJB equation of analytic nonlinear systems becomes an infinite-dimensional Riccati equation.

\end{abstract}

\section{Introduction}

We revisit a 1961 work by Al'brekht \cite{albrekht1961optimal} on the power
series approach to optimal control. He considered an analytic control system
with an analytic cost and showed how to compute the optimal feedback near an
equilibrium one term at a time. The first step is the linear-quadratic
problem, as the quadratic part of the value function solves a Riccati equation.
Each further term is then obtained from a linear equation whose right-hand side
depends only on the terms already known. In 1969, Lukes
\cite{lukes1969optimal} proved that the series converges and that the feedback
it defines is optimal and stabilizing near the equilibrium, as long as the
linearized system is stabilizable and detectable. Over the years, this
recursion has become the standard tool to compute optimal feedback for smooth
nonlinear systems. It has also been turned into algorithms and software
\cite{navasca2007patchy}. In this paper, we turn the question around. Rather than computing the terms of the series individually, we ask what equation they satisfy together. Carleman linearization gives
the answer. After lifting, all the terms of the series together solve a single
infinite-dimensional Riccati-type equation, whose unknown is an infinite matrix.
The first entry of this equation is the Riccati equation of the standard
linearization at the equilibrium.The lifted equation has a structure of its own that characterizes the exact solution, and this paper exploits it.

The series comes from the Hamilton--Jacobi--Bellman (HJB) equation. Optimal
control seeks a control law that minimizes a given cost functional
\cite{khalil2002nonlinear}. Dynamic programming \cite{bellman1954dynamic}
reduces this problem to the HJB equation, a nonlinear partial differential
equation whose solution is the value function \cite{kalman1963theory}. For
nonlinear systems, the equation rarely has a closed-form solution, and a
smooth solution may not exist globally \cite{bardi1997optimal}. One therefore
turns to approximations, such as Galerkin projection \cite{beard1997galerkin}
and polynomial approximation \cite{kalise2018polynomial}. The power series
method of Al'brekht and Lukes is one of them. It computes the Maclaurin
coefficients of the value function one degree at a time. It has been applied
to bilinear systems \cite{cebuhar1984approximation} and to systems with
additive control \cite{yoshida1989quadratic}. The method is local by nature.
Adding terms improves the accuracy near the equilibrium, but it does not
enlarge the region where the series is valid \cite{navasca2007patchy}.

Carleman linearization lifts a nonlinear finite-dimensional system to an
infinite-dimensional bilinear one
\cite{carleman1932application,kowalski1991nonlinear,amini2022finite}. It is
the dynamical systems analogue of the Maclaurin expansion, as the state is
replaced by the sequence of its monomials, and the nonlinear vector field by an infinite state matrix. The idea extends to Fourier bases
\cite{motee2026carleman,chen2025beyond}. It also underlies quantum algorithms for nonlinear differential equations \cite{liu2021efficient,jennings2025quantum}. Related lifts to higher-dimensional linear models have been used to compute
invariant sets \cite{wang2023computation}. Carleman linearization has also
been used to extend the Lyapunov equation to nonlinear systems
\cite{hernandez2004generalization} and to estimate domains of attraction
\cite{loparo1978estimating}. In both cases, the lifted equation is linear in the unknown. The price of the lift is that the resulting operators are
unbounded. In practice, they must be truncated
\cite{forets2017explicit,amini2022finite}.

In earlier work, we lifted the HJB equation itself
\cite{amini2020approximate,amini2020quadratization}. This produced an
infinite-dimensional equation, quadratic and symmetric in the unknown, that
resembles an algebraic Riccati equation. We also gave an iterative algorithm
to solve it. However, one term of that equation contained an operator that
depends on the unknown only implicitly. This operator concealed the structure of the
equation. Since then, truncations of the lifted equation have been learned
from data \cite{kar2023rl,zhao2026data}. A related Riccati-like operator
equation has also been derived in a Koopman framework
\cite{breiten2025unifying}. None of these works gives the solution of the
lifted equation in explicit form. Without it, one cannot check what the
equation actually computes. This is where the scalar case helps. For a scalar
system, the HJB equation is quadratic in the derivative of the value function
and has a closed-form solution. Every claim about the lifted equation can
therefore be checked exactly, as was done for scalar bilinear systems in
\cite{cebuhar1984approximation}. We focus on the scalar case, where the structure of the lifted
equation can be shown explicitly, as a first step towards multivariable systems.

For scalar systems, our contributions are threefold. First, we eliminate the
implicit operator. The lifted HJB equation then becomes a Riccati-type
equation for a lower-triangular Toeplitz matrix. Second, we show that, under
controllability of the linearization and detectability of the cost, this
equation has a unique solution with positive diagonal entries. The solution is
obtained from the scalar Riccati equation of the linearized problem, followed
by a sequence of scalar linear equations. It is a bounded operator whenever
the Maclaurin coefficients of the derivative of the value function decay fast
enough. Third, we prove that the finite sections of this solution converge
exponentially to the derivative of the value function near the origin.
Numerical simulations illustrate the results, and Fig.~\ref{fig:concept} summarizes the approach and the role of the theorems.

\section{Problem Statement}\label{sec:problem}
We start from Al'brekht's 1961 problem \cite{albrekht1961optimal} and
specialize it to one state and one input, which is its simplest nonlinear case.
Its solution has a closed form, which serves as a benchmark as it tells us
what the lifted HJB equation must compute. In this section, we state the
problem, list our assumptions, and give the solution.

The system is input-affine,
\begin{align}
    \dot{x} = f(x) + g(x)u, \qquad x(0)=x_0\in\mathbb{R},
    \label{eqn:GeneralNonlinear}
\end{align}
with input $u\in\mathbb{R}$. We assume that $f(0)=0$ and that $f$ and $g$ are
analytic around the origin, with Maclaurin expansions
\begin{align}
    f(x)=\sum_{k=1}^\infty f_kx^k \ \ {\rm and}\ \  g(x) = \sum_{k=0}^\infty g_kx^k.
    \label{eqn:taylorAssumption}
\end{align}
The cost to be minimized is
\begin{align}
    J(x_0,u)=\int_0^\infty \big(q(x(t)) + ru^2(t)\big)\,dt,
    \label{eqn:costFunctional}
\end{align}
where $r>0$ weights the control and the state cost $q$ satisfies the following
condition.

\medskip
\noindent {\bf Detectability.} The state cost $q:\mathbb{R}\to\mathbb{R}_+$
is nonnegative, and $f$ and $q$ have no common zero other than the origin, i.e., $q(x)\ge 0$ for all $x$ and $f(x)=q(x)=0$ only if $x=0$. In addition, $q$ is analytic around the origin with
\begin{equation}
\label{costfunction.expansion}
q(x)=\sum_{k=2}^\infty q_k x^k  \quad {\rm with} \quad  q_2>0.
\end{equation}
The second part of \eqref{costfunction.expansion}, along with the condition that $f$ and $q$ have no common zero other than the origin, says that the origin is the
only rest point of zero cost. This is zero-state detectability of $(f,q)$
\cite{byrnes1991passivity}. The condition $q_2>0$ means that the linearized
pair $(\sqrt{q_2},f_1)$ is observable.

\medskip
Let $V(x_0)=\inf_u J(x_0,u)$ be the value function, and let $V_x$ denote its derivative. Dynamic programming gives
\begin{equation}\label{eqn:HJBmin}
\min_{u\in\mathbb{R}}\big\{V_x\,\big(f(x)+g(x)u\big)+q(x)+ru^2\big\}=0 .
\end{equation}
The expression in braces is a convex quadratic in $u$. Its minimizer is the
optimal feedback
\begin{align}
    u^*(x)=-\frac{1}{2r}\,g(x)V_x .
    \label{eqn:uStar}
\end{align}
Substituting \eqref{eqn:uStar} into \eqref{eqn:HJBmin} gives the HJB equation
\begin{align}
    V_x f(x)-\frac{1}{4r} (V_x g(x))^2+q(x)=0.
    \label{eqn:nonlinearHJB}
\end{align}
The optimal feedback is known as soon as $V_x$ is. The rest of the paper is
about $V_x$.

\medskip
\noindent {\bf Controllability.} The linearization of
\eqref{eqn:GeneralNonlinear} at the origin is $\dot x=f_1x+g_0u$. We assume
that it is controllable,
\begin{equation}\label{g0.requirement}
g_0\neq 0 .
\end{equation}

\medskip
\noindent {\bf Existence.} Let us define 
\begin{equation} \label{g0f1.requirement}
\Delta:=f_1^2+\frac{1}{r}\,q_2g_0^2 .
\end{equation}
By \eqref{costfunction.expansion} and \eqref{g0.requirement}, $\Delta>0$. The
numbers $\pm\sqrt{\Delta}$ are the eigenvalues of the Hamiltonian matrix of
the linearized problem. The coefficient of $x^2$ in \eqref{eqn:nonlinearHJB}
is its Riccati equation,
\begin{equation}\label{eqn:scalarRiccati}
q_2+f_1v-\frac{g_0^2}{4r}\,v^2=0, \qquad v=V_{xx}(0).
\end{equation}
Its roots are $v=2rg_0^{-2}(f_1\pm\sqrt{\Delta})$. They are real and have
opposite signs, because $q_2>0$. The positive root,
$v=2rg_0^{-2}(f_1+\sqrt{\Delta})$, gives the closed-loop pole
$f_1-\frac{g_0^2}{2r}v=-\sqrt{\Delta}$. As done by Al'brekht, the recursion shown in this paper starts from the mentioned root.

\medskip
\noindent {\bf Closed-form solution.} The recursion is not needed to find
$V_x$. The HJB equation \eqref{eqn:nonlinearHJB} is quadratic in the single
unknown $V_x$, and it can be solved in closed form. Let
\begin{equation}\label{l.def}
l(x):=x^{-2} f(x)^2+r^{-1}x^{-2} q(x)g(x)^2 .
\end{equation}
By \eqref{eqn:taylorAssumption} and \eqref{costfunction.expansion}, $l$ is
analytic around the origin with $l(0)=\Delta>0$. By
\eqref{costfunction.expansion}, $l(x)>0$ for $x\neq0$ wherever $g(x)\neq0$.
Among the solutions of \eqref{eqn:nonlinearHJB} that are twice differentiable
at the origin, exactly one has $V_{xx}(0)>0$:
\begin{align}
    V_x= \frac{2r \big(f(x)+x \sqrt{l(x)}\big)}{ g(x)^2} .
    \label{eqn:V(x)}
\end{align}
The other one has $V_{xx}(0)=2rg_0^{-2}(f_1-\sqrt{\Delta})<0$. To see this,
note that for $x\neq0$ a solution must coincide with one of the two branches
$2r\big(f\pm\sqrt{f^2+qg^2/r}\big)/g^2$, which are distinct because $l>0$. If
$V_{xx}(0)$ exists, the one-sided slopes at the origin must agree, and this leaves the two analytic functions $2r(f\pm x\sqrt{l})/g^2$. By
\eqref{eqn:uStar}, the optimal feedback is 
\begin{align}
    u^*(x)&=-\frac{f(x)+x\sqrt{l(x)}}{g(x)} \nonumber
\end{align}
and the optimal closed-loop system is 
\begin{align}
\label{closedloop.explicit}
\dot{x} &= f(x)+g(x)u^*(x)=-x\sqrt{l(x)}.
\end{align}

\medskip
On any interval around the origin where $g\neq0$, the closed loop in
\eqref{closedloop.explicit} strictly decreases $|x|$. The origin is therefore
asymptotically stable under $u^*$, with linearization
$\dot x=-\sqrt{\Delta}\,x$. The function $V(x)=\int_0^xV_x(s)\,ds$ is positive
definite. Since $f$, $g$ and $l$ are analytic around the origin with
$g(0)=g_0\neq0$ and $l(0)=\Delta>0$, the function $V_x$ in \eqref{eqn:V(x)} is
analytic around the origin:
\begin{equation}\label{Vx.expansion}
V_x=\sum_{k=1}^\infty v_k x^k, \qquad x\in (-\rho, \rho).
\end{equation}
Here $\rho$ is any number below the radius of convergence. By Cauchy's
estimate, there is a constant $M=M(\rho)$ such that
\begin{equation}
|v_k|\le M \rho^{-k},\qquad k\ge 1 .
\label{vk.def}
\end{equation}
The radius of convergence is set by the complex singularities of $V_x$, such
as the complex zeros of $g$ and of $l$. It can be smaller than the real
interval on which \eqref{eqn:V(x)} is defined. We return to this point in
Section~\ref{sec:finitesection}.

While the closed form solves the scalar problem, it does not show how the
coefficients $v_k$ fit together. In the next section, we lift the HJB equation
by Carleman linearization. We then show that all of the coefficients $v_k$ solve a single
infinite-dimensional Riccati-type equation whose first entry is the Riccati equation
\eqref{eqn:scalarRiccati}.

\medskip

\section{HJB Equation for Carleman Linearization}
\label{HJBCarleman.section}
\begin{figure}[h]
    \centering
    \includegraphics[width=\linewidth,]{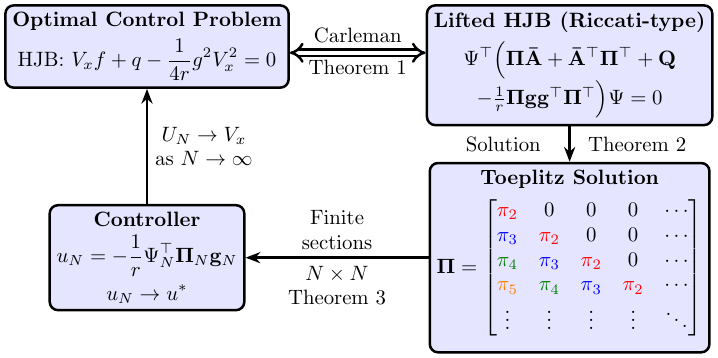}
    \caption{The flowchart of our approach: (i) The Carleman lift $\Psi(x)$ maps the \emph{Optimal Control Problem} to the \emph{Lifted HJB} equation, whose unknown is a lower-triangular Toeplitz matrix (Theorem~\ref{HJBreformulation.theorem} in Section \ref{HJBCarleman.section}). (ii) The \emph{Toeplitz solution} follows from one scalar Riccati equation and then a sequence of  linear equations (Theorem~\ref{HJBcoefficient.thm} in Section \ref{sec:finitesection}). (iii) Finite sections of the Toeplitz solution give the \emph{Controller}, converging exponentially near the origin (Theorem~\ref{finitesection.thm} in Section \ref{sec:finitesection}). 
    }
\label{fig:concept}
\vspace{-0.4cm}
\end{figure}

Carleman linearization replaces the scalar state $x$ by the sequence of its
powers \cite{carleman1932application,kowalski1991nonlinear}. In the new
variables, the nonlinear system \eqref{eqn:GeneralNonlinear} becomes an
infinite-dimensional bilinear system. We first write this system down and
then lift the HJB equation along with it. The Carleman lifting operator is the infinite-dimensional nonlinear map
\begin{equation}
\label{Carlemanstate.def}
\Psi(x):=\text{col}[x,x^2,x^3,\dots],
\end{equation}
where $\text{col}[v,w]$ stacks the column vector $v$ on top of the vector $w$. Along a trajectory of \eqref{eqn:GeneralNonlinear}, each power satisfies
$\frac{d}{dt}x^m=m\,x^{m-1}\big(f(x)+g(x)u\big)$. Expanding $x^{m-1}f(x)$ and
$x^{m-1}g(x)$ by \eqref{eqn:taylorAssumption} gives the Carleman
linearization
\begin{align}
    \dot{\Psi}={\bf D} \big(\mathbf{A}\Psi + (\mathbf{B}\Psi + \mathbf{b}_0)u\big),
    \label{eqn:LiftedBilinearized}
\end{align}
where  
the operators are defined by
\begin{align}
\mathbf{A} = \begin{bmatrix}
f_1 & f_2 & f_3 & \cdots\\
0 & f_1 & f_2 & \cdots\\
0 & 0 & f_1 & \cdots\\
\vdots & \vdots & \vdots & \ddots
\end{bmatrix},~\mathbf{B} = \begin{bmatrix}
g_1 & g_2 & g_3 & \cdots\\
g_0 & g_1 & g_2 & \cdots\\
0 & g_0 & g_1 & \cdots\\
\vdots & \vdots & \vdots & \ddots
\end{bmatrix} \label{A.def}
\end{align}
\begin{align}
\mathbf{D} = \begin{bmatrix}
1 & 0 & 0 & \cdots\\
0 & 2 & 0 & \cdots\\
0 & 0 & 3 & \cdots\\
\vdots & \vdots & \vdots & \ddots
\end{bmatrix}, \qquad
\mathbf{b}_0 = g_0\mathbf{e}_1 = \begin{bmatrix} g_0\\ 0\\ 0\\ \vdots \end{bmatrix} \label{B.def}
\end{align}
and $\mathbf{e}_1=[1,0,0,\dots]^\top$.
The structure of these matrices is the
point of this section. $\mathbf{A}$ is upper-triangular Toeplitz, with the
Maclaurin coefficients of $f$ along its rows. $\mathbf{B}$ is upper Hessenberg
Toeplitz, with the coefficients of $g$ along its rows and $g_0$ on the
subdiagonal. On the lifted state they act as multiplication operators:
\begin{equation}\label{APsi.BPsi.eqn}
\mathbf{A}\Psi = x^{-1}f(x)\,\Psi,
\qquad
\mathbf{B}\Psi+\mathbf{b}_0 = x^{-1}g(x)\,\Psi ,
\end{equation}
because the $m$-th entries of the two sides are $x^{m-1}f(x)$ and
$x^{m-1}g(x)$. These two identities are all we need from the lift. In this
section, we study the HJB equation of \eqref{eqn:LiftedBilinearized}, and in
 Section~\ref{sec:finitesection} we prove exponential convergence of its
finite-section approximations; see Figure~\ref{fig:concept}.

The cost lifts in the same way. By \eqref{costfunction.expansion},
\begin{equation}
\label{Q.def}
q(x)=  \Psi^\top {\bf Q} \Psi,
\end{equation}
where ${\bf Q}=[q(i,j)]_{i,j=1}^\infty$ is the symmetric tridiagonal matrix
with $q(i,i)=q_{2i}$, $q(i,i+1)=q(i+1,i)=q_{2i+1}/2$ for $i\ge1$, and
$q(i,j)=0$ otherwise. Following the quadratization procedure in
\cite{amini2020quadratization}, we view the value function as a function of
$\Psi$ and write $V_\Psi$ for its gradient. The HJB equation for the lifted
system \eqref{eqn:LiftedBilinearized} is then
\begin{align}
   V_{\Psi}^\top {\bf D}\mathbf{A}\Psi
   - \frac{1}{4r}\big(V_{\Psi}^\top {\bf D}(\mathbf{B}\Psi+\mathbf{b}_0)\big)^2
   + \Psi^\top \mathbf{Q} \Psi=0,
 \label{HJB.linearization}
\end{align}
and the corresponding optimal control law is
\begin{align}
    u^* = -\frac{1}{2r}(\mathbf{B}\Psi+\mathbf{b}_0)^\top {\bf D}V_{\Psi}.
    \label{eqn:bilinearHJB}
\end{align}
By the chain rule, $V_x=V_\Psi^\top\,\frac{d\Psi}{dx}=x^{-1}V_\Psi^\top{\bf D}\Psi$.
Together with \eqref{APsi.BPsi.eqn}, this shows that on the lifted state
\eqref{HJB.linearization} and \eqref{eqn:bilinearHJB} are the HJB equation
\eqref{eqn:nonlinearHJB} and the feedback \eqref{eqn:uStar} of
Section~\ref{sec:problem}. Lifting therefore leaves the HJB equation unchanged on the lifted state and expresses it through the lifted operators, whose structure is exploited next.

We look for a solution whose gradient is linear in $\Psi$,
\begin{align}
    {\bf D}V_\Psi=2{\bf P} \Psi
    \label{P.def}
\end{align}
for some infinite matrix ${\bf P}=[p(i,j)]_{i,j=1}^\infty$. Then, 
$V_x=2x^{-1}\Psi^\top{\bf P}\Psi$. Since
$\Psi^\top{\bf P}\Psi=\sum_{i,j=1}^\infty p(i,j)x^{i+j}$, the matrix ${\bf P}$ enters
only through the sums of its anti-diagonal entries,
\begin{equation}\label{pim.def}
\pi_m({\bf P})=\pi_m({\bf P}^\top)=\sum_{i=1}^{m-1} p(i, m-i)
\end{equation}
for all $m\ge 2$ and $V_x=2\sum_{m=2}^\infty\pi_m({\bf P})\,x^{m-1}$. Comparing with
\eqref{Vx.expansion}, the Maclaurin coefficients of $V_x$ are
$v_{m-1}=2\pi_m({\bf P})$. We collect the anti-diagonal sums in the
lower-triangular Toeplitz matrix
\begin{align}
    \mathbf{\Pi}:=\begin{bmatrix}
        \pi_2({\bf P}) &0&0&\cdots\\
        \pi_3 ({\bf P}) & \pi_2 ({\bf P}) &0 & \cdots\\
        \pi_4 ({\bf P})& \pi_3 ({\bf P}) & \pi_2 ({\bf P}) & \cdots\\
        \vdots & \vdots & \vdots & \ddots
    \end{bmatrix}.
    \label{Pi.def}
\end{align}
Many matrices ${\bf P}$ have the same anti-diagonal sums. They all give the
same $V_x$, and $\mathbf{\Pi}$ is what they have in common. The next theorem
rewrites the HJB equation \eqref{HJB.linearization} and the control law
\eqref{eqn:bilinearHJB} in terms of $\mathbf{\Pi}$ alone. Let us denote 
${\bf f}=[f_1,f_2,f_3,\ldots]^\top$, ${\bf g}=[g_0,g_1,g_2,\ldots]^\top$, and
$\bar{\bf A}={\bf f}\,{\bf e}_1^\top$.

\vspace{0.2cm}
\begin{theorem} \label{HJBreformulation.theorem}
Suppose that ${\bf A}$, ${\bf B}$, ${\bf P}$, $\mathbf{Q}$, and $\mathbf{\Pi}$ are as in \eqref{A.def},
\eqref{B.def}, \eqref{P.def}, \eqref{Q.def}, and \eqref{Pi.def}. Under the ansatz \eqref{P.def}, the
lifted HJB equation \eqref{HJB.linearization} is equivalent to
\begin{align}
    \Psi^\top
    \Big( \mathbf{\Pi}\bar{\bf A}+\bar{\bf A}^\top \mathbf{\Pi}^\top + {\bf Q}
    - \frac{1}{r}\mathbf{\Pi}\mathbf{g}\mathbf{g}^\top\mathbf{\Pi}^\top\Big)\Psi=0,
    \label{HJB.linearization.invariant}
\end{align}
and the control law \eqref{eqn:bilinearHJB} becomes
\begin{equation} \label{eqn:bilinearHJB.invariant}
u^*=-\frac{1}{r}\, \Psi^\top \mathbf{\Pi}\, {\bf g} .
\end{equation}
\end{theorem}
\vspace{0.2cm}

The proof is given in Appendix~\ref{app:proofThm1}. Before we use the
theorem, we point out what it implies. Equation \eqref{HJB.linearization.invariant} has the form of an algebraic
Riccati equation, $PA+A^\top P+Q-\frac1r Pbb^\top P=0$, with $\mathbf{\Pi}$
in the role of $P$. The data are read off the lift where ${\bf f}$ is the first
row of ${\bf A}$, ${\bf g}$ is the second row of ${\bf B}$, and ${\bf Q}$ is
the lifted cost. The unknown is a single sequence, $\pi_2,\pi_3,\dots$, which, for the solution of Theorem~\ref{HJBcoefficient.thm}, is the sequence of Maclaurin coefficients of $V_x$ up to the factor
two. In our earlier work \cite{amini2020quadratization}, the quadratic term of
the lifted equation contained an auxiliary operator that depended on the
unknown only implicitly. Here that operator is gone. The quadratic term is
$\mathbf{\Pi}\mathbf{g}\mathbf{g}^\top\mathbf{\Pi}^\top$, explicit and of
rank one, as in the single-input linear-quadratic case; along $\Psi(x)$ it
gives $\frac1r(\Psi^\top\mathbf{\Pi}\mathbf{g})^2=r\,(u^*)^2$, the control cost.
The $(1,1)$ entry of the matrix in \eqref{HJB.linearization.invariant} is
$2f_1\pi_2+q_2-\frac{g_0^2}{r}\pi_2^2$. With $v=2\pi_2$, this is the Riccati
equation \eqref{eqn:scalarRiccati} of the linearization. The infinite Riccati
equation thus contains the finite one as its first entry, and the rest of its entries hold the nonlinear terms.

Two facts explain why the lift produces an equation with a Riccati form rather than
something more complicated. First, by \eqref{APsi.BPsi.eqn}, the Carleman
vector $\Psi(x)$ is an eigenvector of ${\bf A}$ with eigenvalue $f(x)/x$, and
the analogous relation holds for ${\bf B}$. Thus, the lifted equation depends on ${\bf P}$ only through the series
$\sum_{m=2}^\infty\pi_m({\bf P})x^m$. Second, products of power series are
represented by products of lower-triangular Toeplitz matrices. The matrix
$\mathbf{\Pi}$ is the Toeplitz matrix of the series
$\sum_{m=2}^\infty\pi_m({\bf P})x^{m-2}$, and for any sequence
${\bf c}=[c_0,c_1,\dots]^\top$ with generating function $c(x)$,
\begin{equation}\label{Pic.identity}
\Psi^\top\mathbf{\Pi}\,{\bf c}=x\left(\sum_{m=2}^\infty\pi_m({\bf P})\,x^{m-2}\right)c(x).
\end{equation}
This identity turns every product in the HJB equation into a matrix product;
both identities in the proof of Theorem~\ref{HJBreformulation.theorem} are
instances of it. We should highlight that \eqref{HJB.linearization.invariant} is an identity between quadratic
forms evaluated along the curve $\Psi(x)$. It says that every anti-diagonal
sum of the matrix in parentheses vanishes, which is one scalar equation per
degree. It does not say that the matrix vanishes entry by entry.

Theorem~\ref{HJBreformulation.theorem} writes the lifted HJB equation in
terms of $\mathbf{\Pi}$ alone, but \eqref{HJB.linearization.invariant} is
still infinite-dimensional. The next corollary closes the loop with
Section~\ref{sec:problem}. The first column of $\mathbf{\Pi}$ defines a
function that solves the scalar HJB equation \eqref{eqn:nonlinearHJB}, and
the lifted feedback \eqref{eqn:bilinearHJB.invariant} is the feedback
\eqref{eqn:uStar}. This is what allows us to compute $\mathbf{\Pi}$ degree by
degree in Section~\ref{sec:finitesection} and to compare its finite sections
with $V_x$ in \eqref{eqn:V(x)}.

\vspace{0.2cm}
\begin{corollary}\label{U.cor}
Let $\mathbf{\Pi}$ be a lower-triangular Toeplitz matrix of the form
\eqref{Pi.def} that satisfies \eqref{HJB.linearization.invariant}, and set
\begin{equation}\label{U.def}
U(x) := 2\,\Psi^\top\mathbf{\Pi}\,\mathbf{e}_1
= 2\sum_{m=2}^\infty \pi_m(\mathbf{P})\,x^{m-1} .
\end{equation}
Then, $U(x)=2x^{-1}\Psi^\top{\bf P}\Psi$ and the function $U$ satisfies the HJB
equation
\begin{equation}
\label{U.HJB}
q(x)+U(x) f(x)-\frac{1}{4r} \big(g(x) U(x)\big)^2=0,
\end{equation}
and the control \eqref{eqn:bilinearHJB.invariant} equals
$u^*=-\frac{1}{2r}g(x)U(x)$, as in \eqref{eqn:uStar}. Conversely, let $U$ be
analytic around the origin with $U(0)=0$ and satisfy \eqref{U.HJB}. Then, the
Toeplitz matrix $\mathbf{\Pi}$ whose entries $\pi_m$ are one half of the
coefficient of $x^{m-1}$ in $U$ satisfies \eqref{HJB.linearization.invariant}.
\end{corollary}
\vspace{0.2cm}

The proof is given in Appendix~\ref{app:proofCor1}. In words, the lifted Riccati-type equation and the scalar HJB equation have the same solutions.

\vspace{0.2cm}
\begin{remark}
By \eqref{P.def} and the chain rule,
\begin{align}
\label{VxU.eqn}
    V_x=\frac{dV}{d\Psi}\frac{d\Psi}{dx}= x^{-1} V_\Psi^\top\mathbf{D} \Psi= 2x^{-1} \Psi^\top {\bf P} \Psi =U(x).
\end{align}
Thus, the function $U$ in \eqref{U.def} is the derivative of the lifted value
function along $\Psi(x)$. Along the closed loop, \eqref{U.HJB} and
\eqref{eqn:bilinearHJB.invariant} give
\begin{align}
    \dot{V}&= V_x \big(f(x)+g(x) u^*\big)
    = -\Psi^\top\mathbf{Q}\Psi-\frac{1}{r}\,\Psi^\top\mathbf{\Pi}\mathbf{g}\mathbf{g}^\top\mathbf{\Pi}^\top\Psi
    \nonumber\\
    &=-q(x)-r\,|u^*(x)|^2 \le 0 .
    \nonumber
\end{align}
This is the Lyapunov inequality one expects from an HJB equation. Moreover, the solution constructed in Theorem~\ref{HJBcoefficient.thm} has
$\pi_2({\bf P})>0$ by \eqref{HJBcoefficient.thm.eq1}, and $U$ is analytic
with $U'(0)=2\pi_2({\bf P})>0$. By Section~\ref{sec:problem}, there is only
one such solution of the HJB equation. Hence, $U$ is the closed-form $V_x$ in
\eqref{eqn:V(x)}, the closed-loop system is $\dot x=-x\sqrt{l(x)}$ as in
\eqref{closedloop.explicit}, the origin is asymptotically stable in the sense of Lyapunov, and $V$ is
the value function of the problem. Thus, the positive-diagonal solution of the lifted Riccati-type equation gives the optimal cost while serving as a Lyapunov function.
\end{remark}

\section{Exponential Convergence of Finite-Section Approximations}
\label{sec:finitesection}

By Theorem~\ref{HJBreformulation.theorem}, solving the lifted HJB equation
\eqref{HJB.linearization} amounts to finding the lower-triangular Toeplitz
solution $\mathbf{\Pi}$ of the Riccati-type equation
\eqref{HJB.linearization.invariant}. In this section, we solve this equation
by a recursion of scalar equations, and we show that the finite-sections of
the solution converge exponentially to $V_x$ near the origin.

\vspace{0.2cm}
\begin{theorem}
\label{HJBcoefficient.thm}
Suppose that the system \eqref{eqn:GeneralNonlinear} with $f$, $g$ and $q$
satisfying \eqref{eqn:taylorAssumption}, \eqref{costfunction.expansion} and
\eqref{g0.requirement}, and let $\Delta>0$ be as in \eqref{g0f1.requirement}.
Then, the Riccati-type equation \eqref{HJB.linearization.invariant} has
exactly one lower-triangular Toeplitz solution $\mathbf{\Pi}$ with positive
diagonal entries, where its diagonal entry is
\begin{equation}
 \label{HJBcoefficient.thm.eq1}
\pi_2({\bf P})=
r g_0^{-2} (f_1+\sqrt{\Delta}),
\end{equation}
and its remaining entries $\pi_{m+2}({\bf P})$, for $m\ge 1$, are given
recursively by
\begin{align}   \label{HJBcoefficient.thm.eq2}
 & \hspace{-0.2cm} \pi_{m+2}({\bf P}) =   (2\sqrt{\Delta})^{-1}\Big( q_{m+2}\nonumber \\
& \hspace{-0.2cm} \ \ +
2\sum_{k=0}^{m-1} \big(f_{m-k+1}- g_0 g_{m-k}\pi_2({\bf P})/r\big) \pi_{k+2}({\bf P})\nonumber\\
& \hspace{-0.2cm} \ \ -
\frac{1}{r} \sum_{n=1}^{m-1}
       \sum_{j=0}^n \sum_{j'=0}^{m-n} g_{n-j} g_{m-n-j'}
\pi_{j+2}({\bf P}) \pi_{j'+2}({\bf P})\Big).
\end{align}
Moreover, $2\pi_m({\bf P})=v_{m-1}$ for all $m\ge2$, where $v_k$ are the
Maclaurin coefficients of $V_x$ in \eqref{Vx.expansion}.
\end{theorem}

\vspace{0.2cm}
The proof is given in Appendix~\ref{app:proofThm2}. The recursion deserves
a closer look. Its first equation, the case $m=0$, is the Riccati equation
\eqref{eqn:scalarRiccati} of the linearization, and
\eqref{HJBcoefficient.thm.eq1} is its positive  solution. 
Every later equation is
linear in the new unknown $\pi_{m+2}({\bf P})$, and the coefficient in front
of it is always the same number, $-2\sqrt{\Delta}$, which is twice the
closed-loop pole of Section~\ref{sec:problem}. The right-hand side involves
only the data $f_k$, $g_k$, $q_k$ and the entries already computed. This is
Al'brekht's recursion for the scalar case, with one difference: {\it here it is
read off a single matrix equation, degree by degree, rather than derived
anew at each degree.} Two consequences follow. Existence and uniqueness of
the solution are decided entirely by the first equation, since $\Delta>0$
makes every later pivot nonzero. And, the first $N$ equations involve only
$\pi_2({\bf P}),\dots,\pi_{N+1}({\bf P})$. Therefore, truncating the equation and
truncating the solution give the same result. The finite-section of the
solution is the exact solution of the finite-section. No truncation error is
introduced by the recursion.

Let us consider the finite-section of order $N\ge 1$,
\begin{align}
\label{PiN.def}
\mathbf{\Pi}_N=\scalebox{0.9}{$\begin{bmatrix}
        \pi_2({\bf P}) &0&0&\cdots & 0\\
        \pi_3 ({\bf P}) & \pi_2 ({\bf P}) &0 & \cdots & 0 \\
        \pi_4 ({\bf P})& \pi_3 ({\bf P}) & \pi_2 ({\bf P}) & \cdots & 0\\
        \vdots & \vdots & \vdots & \ddots & \ddots\\
         \pi_{N+1} ({\bf P})& \pi_{N} ({\bf P}) & \pi_{N-1} ({\bf P}) & \cdots & \pi_2({\bf P})
    \end{bmatrix}$},
\end{align}
of the matrix $\mathbf{\Pi}$ in \eqref{Pi.def}, and define
\begin{equation}
\label{UN.def}
U_N(x)= 2\Psi_N^\top\mathbf{\Pi}_N\mathbf{e}_1,
\end{equation}
where $\Psi_N=[x, x^2, \ldots, x^N]^\top$. By Theorem~\ref{HJBcoefficient.thm},
$U_N(x)=\sum_{k=1}^{N}v_kx^k$ is the Taylor polynomial of degree $N$ of $V_x$.
The next theorem gives the rate at which $U_N$ converges to $V_x$.

\vspace{0.2cm}
\begin{theorem}\label{finitesection.thm}
Suppose that the system \eqref{eqn:GeneralNonlinear} with $f$, $g$ and $q$
satisfying \eqref{eqn:taylorAssumption}, \eqref{costfunction.expansion} and
\eqref{g0.requirement}, and let $V_x$ be the solution \eqref{eqn:V(x)} of the
HJB equation \eqref{eqn:nonlinearHJB}. Let $\pi_m({\bf P})$, for $m\ge 2$, be as
in \eqref{HJBcoefficient.thm.eq1} and \eqref{HJBcoefficient.thm.eq2}, and
let $\mathbf{\Pi}_N$ and $U_N$ be as in \eqref{PiN.def} and \eqref{UN.def}.
Then, 
\begin{align}\label{finitesection.thm.eq1}
|U_N(x)-V_x|\le M (1-|x|/\rho)^{-1} |x/\rho|^{N} 
\end{align}
for all $x\in (-\rho, \rho)$ and $N\ge 1$, where $M$ and $\rho$ are the
constants in \eqref{vk.def}.
\end{theorem}

The proof is given in Appendix~\ref{app:proofThm3}. The bound is exponential
in $N$ with ratio $|x|/\rho$. Each additional row of $\mathbf{\Pi}_N$ reduces
the error by this factor. The constant $\rho$ is any number below the radius
of convergence of $V_x$, and, as noted in Section~\ref{sec:problem}, that
radius is set by the complex singularities of $V_x$. It can be smaller than
the real interval on which \eqref{eqn:V(x)} is defined, in which case the
finite-sections converge only on $(-\rho,\rho)$ even though $V_x$ is defined
beyond it. The truncated feedback $u_N^*(x):=-\frac{1}{r}g(x)U_N(x)$
inherits the bound with an extra factor $|g(x)|/(2r)$.

The state matrix ${\bf D}{\bf A}$ arising in the Carleman linearization
\eqref{eqn:LiftedBilinearized} is unbounded, whereas the 
Riccati-type solution $\mathbf{\Pi}$ of the corresponding HJB equation need not be. In fact,
as shown in the following remark, $\mathbf{\Pi}$
is  a bounded operator whenever the Maclaurin coefficients of $V_x$
exhibit sufficiently fast decay. This distinction highlights that the
unboundedness of the lifted dynamics does not preclude the existence of a
well-behaved solution to the associated 
 HJB equation.

\vspace{0.2cm}

\begin{remark}
For an infinite matrix ${\bf C}=[c(i,j)]_{i,j=1}^\infty$, define its Schur
norm by
$$\|{\bf C}\|_{\mathcal S}=\max \Bigg\{\sup_{i\ge 1}\sum_{j=1}^\infty |c(i,j)|,\ \sup_{j\ge 1}\sum_{i=1}^\infty |c(i,j)|\Bigg\}. $$
Suppose that the radius of convergence of $V_x$ in \eqref{Vx.expansion}
exceeds one, and choose $\rho>1$ in \eqref{vk.def}. Since $\mathbf{\Pi}$ is
lower-triangular Toeplitz, every column sum equals
$\sum_{m=2}^\infty|\pi_m({\bf P})|$, and by $2\pi_m({\bf P})=v_{m-1}$ and
\eqref{vk.def},
\begin{align}
    \|{\bf\Pi}\|_{\mathcal S}=\sum_{m=2}^\infty |\pi_m({\bf P})|\le \frac{M}{2}
    \sum_{m=2}^\infty \rho^{-(m-1)} = \frac{M}{2(\rho-1)} . \nonumber
\end{align}
By the Schur test, $\mathbf{\Pi}$ is then a bounded operator on $\ell^p$ for
every $1\le p\le \infty$, where $\ell^p$ is the Banach space of $p$-summable
sequences $\{c(k)\}_{k=1}^\infty$ with its standard norm. The threshold one
reflects the choice of unit weights. 
On the weighted space with norm
$\|c\|_{p, \sigma}:=\big(\sum_{k=1}^\infty |c(k)|^p\sigma^{kp}\big)^{1/p}$ for $1\le p<\infty$ and $\|c\|_{\infty, \sigma}=\sup_{k\ge 1} |c(k)|\sigma^k$, the same computation shows that $\mathbf{\Pi}$ is
bounded for every $0<\sigma<\rho$.
\end{remark}

\section{Simulations}
In the previous sections, we showed that the solution of the lifted HJB equation is obtained from one Riccati equation followed by a sequence of linear scalar equations. We also showed that the finite sections $U_N$  for $N\ge 1$ in \eqref{UN.def} converge exponentially to the solution $V_x$ of the HJB equation for $|x|<\rho$, where the convergence radius $\rho$ is related to the Maclaurin series of $V_x$ in Theorem \ref{finitesection.thm}. To show the accuracy of the finite-section approximation, we compare the truncated value gradient and feedback of Theorem \ref{HJBcoefficient.thm},
\begin{align}
    U_N(x)&=2\sum_{i=1}^N \pi_{i+1}
    (\mathbf{P})x^{i}
    \label{eqn:V_N},\hspace{5pt} u_N=-\frac{1}{r}{\Psi}^\top_N\mathbf{\Pi}_N\mathbf{g}_N
\end{align}
against the closed-form solution \eqref{eqn:V(x)} where the truncated vector $\mathbf{g}_N=[g_0,\dots,g_{N-1}]^\top$. We denote by $J_N$ the cost in \eqref{eqn:costFunctional} along the closed loop under $u_N$. We consider three nonlinear systems, where the algorithm in Theorem \ref{HJBcoefficient.thm} was written as a MATLAB program to solve their lifted HJB equation.

\subsection*{Example 1: Polynomial drift with rational input}
\begin{figure}[h]
    \centering
    \includegraphics[width=0.49\linewidth, keepaspectratio]{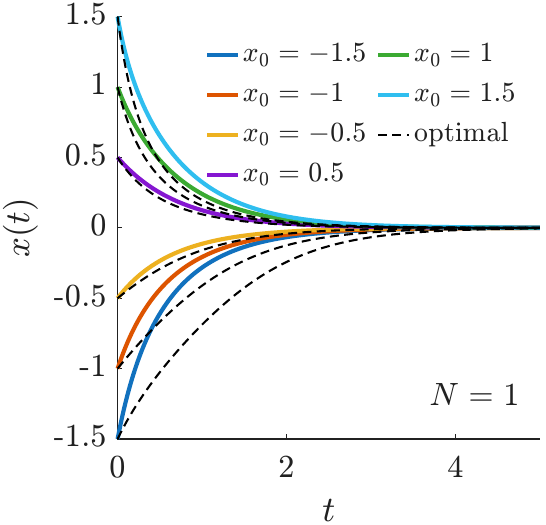}
    \includegraphics[width=0.49\linewidth, keepaspectratio]{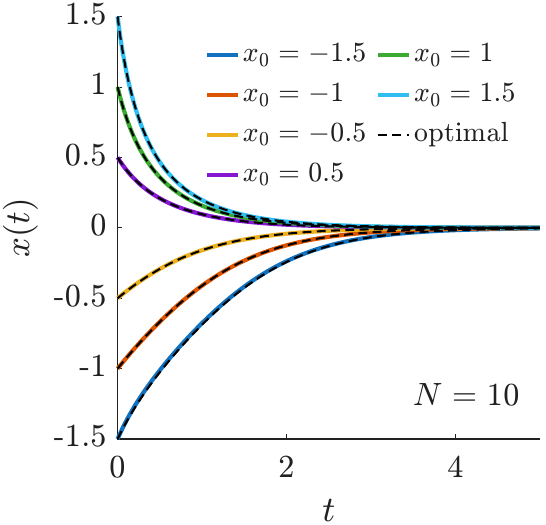}
    \caption{Trajectories of the closed-loop system for Example 1 using $N=1$ (left)
and $N=10$ (right) for different initial conditions.
}
    \label{fig:example_1}
\end{figure}

Consider the nonlinear system
\begin{align}
    \dot{x}=x+x^2+u(1-(x/3))^{-1}
    \label{eqn:example1}
\end{align}
with weight $r=1$ and state cost function $q(x)=x^2$. The drift function $f(x) = x + x^2$ and the control function $g(x) = (1 - x/3)^{-1}$ satisfy \eqref{g0.requirement} and \eqref{g0f1.requirement}. The gradient of the value function and the optimal control law are given by
$$V_x= \frac{2(3x-x^2)}{9}\Big( 3+2x-x^2+ \sqrt{9+(3+2x-x^2)^2}\Big)$$
and
$u^*= -\frac{x}{3}( 3+2x-x^2+ \sqrt{9+(3+2x-x^2)^2})
$
respectively, and the corresponding dynamical system \eqref{eqn:GeneralNonlinear} with optimal control becomes
\begin{equation}
\label{Dynamical.Ex1}
\dot x= -\frac{x  \sqrt{9+(3+2x-x^2)^2}}{3-x}
\nonumber
\end{equation}
In Fig. \ref{fig:example_1}, the truncation order $N$ is held constant, and the closed-loop trajectories are simulated using the truncated controller $u_N$. The left panel in Fig. \ref{fig:example_1} shows the Linear Quadratic Regulator (LQR) approximation, while the improvement is visible in the right panel, where the truncation order is increased to $N=10$ and more coefficients $\pi_m$ enter $\mathbf{\Pi}_N$. It is observed that the approximation yields accurate estimates across different initial conditions for large truncation order $N$. We remark that $V_x$ is analytic on a complex disc $|z|< \sqrt{6-3\sqrt{2}}\approx 1.326$.

\subsection*{Example 2: Trigonometric drift and input}
\begin{figure}[h]
    \centering
    \includegraphics[width=0.47\linewidth]{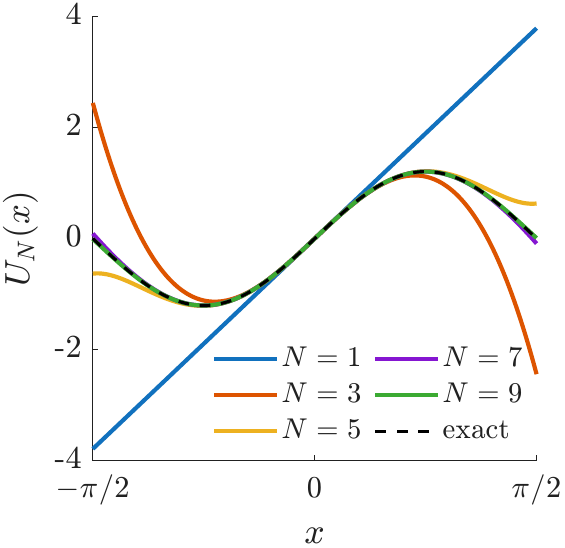}
    \includegraphics[width=0.49\linewidth]{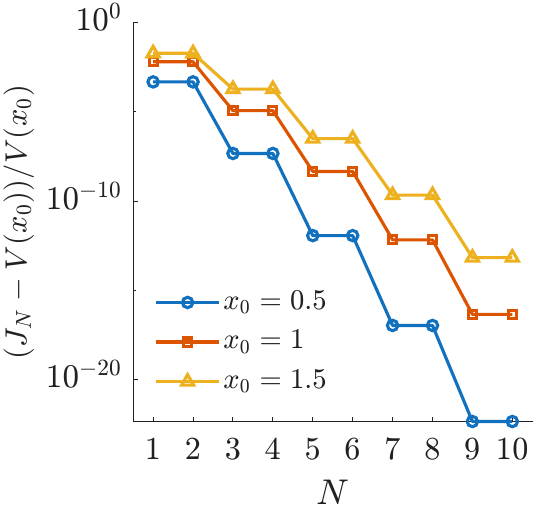}
    \caption{Truncated solution to the HJB compared to the exact solution represented by the dashed line, for different $N$ (left). Log-scale normalized closed-loop cost error against $N$ for different initial conditions (right).}
    \label{fig:example_2}
\end{figure}

Consider the nonlinear system studied by \cite{navasca2007patchy} 
\begin{align}
    \dot{x}=\tan(x)+(\sec(x))u
    \nonumber
\end{align}
with $q(x)=(1/2)\sin^2(x)$ and $r=1/2$, for which the closed-form value function and optimal control are
\begin{align}
    V=\frac{1+\sqrt{2}}{2}\sin^2(x), \quad u^*=-\big(1+\sqrt{2}\big)\sin(x)
    \nonumber
\end{align}
In this case, $f_1=1$, $g_0=1$, and $\Delta=2>0$. Thus, by Theorem \ref{HJBcoefficient.thm}, we have $\pi_2(\mathbf{P})=(1+\sqrt{2})/2$. While $f,g,$ and $q$ are analytic on $|x|<\pi/2$, the poles of $\tan(x)$ and $\sec(x)$ at $\pm \pi/2$ cancel against the factor $2r/g^2=\cos^2(x)$ in \eqref{eqn:V(x)} for this example, giving $V_x=\left(1+\sqrt{2}\right)\sin(x)\cos(x)$, which is analytic everywhere. To verify the approximation, Fig. \ref{fig:example_2} (left) shows that the truncated gradient $U_N(x)$ represents $V_x(x)$ more accurately over the state space as the truncation order is increased. In Fig. \ref{fig:example_2} (right), it is observed that the normalized error between the closed-loop cost of the truncated controller and the optimal cost decreases with more truncations for the listed initial conditions. We note that, since $f$ is odd and both $g$ and $q$ are even, $V_x$ is an odd function, and $\pi_m(\mathbf{P})=0$ for odd $m$. Thus, even truncations equal the preceding odd ones in Fig. \ref{fig:example_2} (right).

\subsection*{Example 3: Fast-intermediate glycolysis model with PFK enzyme activity as input}

\begin{figure}[h]
    \centering
    \includegraphics[width=0.49\linewidth, ]{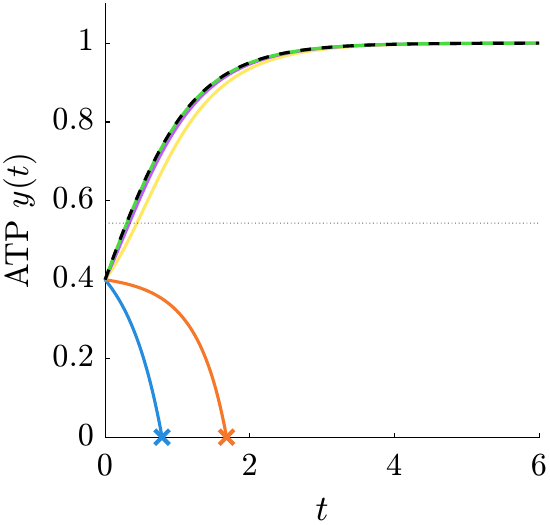}\includegraphics[width=0.49\linewidth, ]{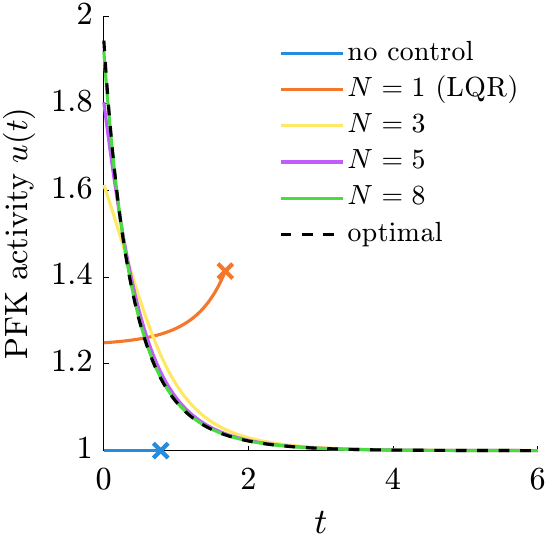}
    \caption{ATP $y(t)$ (left) and PFK activity $u(t)$ (right) from $y(0) =0.40$. Parameters used were $a=1, h=2, c=1, w=1, r=1$. The dashed line is the optimal closed loop, the dotted line in the left panel is the unstable equilibrium $y_u$, and $\times$ marks where the ATP level reaches zero.}
    \label{fig:example_3}
\end{figure}

Consider the autocatalytic glycolysis model studied by \cite{chandra2011glycolytic} where the dynamics are given by 
\begin{align*}
    \dot y = u\,\gamma(y) - c,\qquad \gamma(y)=\frac{2y^{a}}{1+y^{2h}}
\end{align*}
Here, $u$ is the PFK activity, $\gamma$ is the ATP production per unit activity, and $c>0$ is the ATP consumption rate, with $\gamma(1)=1$. With $x=y-1$ and $u=c+v$, this is the general form \eqref{eqn:GeneralNonlinear} with $f=c(\gamma-1), g=\gamma$, and $q=wx^2$, where $w>0$. With $J=\int_0^\infty\big[w\,(y-1)^2+r\,v^2\big]\,dt$, we have $g_0=\gamma(1)=1$ and $\Delta = c^2(a-h)^2+w/r>0$, so the conditions in \eqref{g0.requirement} and \eqref{g0f1.requirement} hold and the HJB optimal controls have the following closed-form
 \begin{align*}
     u^*(y)=\frac{c-\sgn(y-1)S(y)}{\gamma(y)}, \hspace{12pt} S=\sqrt{f^2+\frac{w}{r}\gamma^2(y-1)^2}
 \end{align*}
Thus, $\dot{y}=-\sgn{(y-1)}S(y)$ with $S>0$ for $y\neq1$ and $y>0$ such that the optimal control can recover ATP from any $y(0)>0$. Fig. \ref{fig:example_3} shows ATP $y(t)$ (left) and PFK activity $u(t)$ (right) for the closed-loop system. Next, we discuss the results analytically. ATP rises only when the PFK activity exceeds $c/\gamma(y)$, the level at which ATP production equals consumption. In the uncontrolled case $(u=c)$, we observe an unstable equilibrium $y_u<1$ with $\gamma(y_u)=1$. The optimal control adds $S/\gamma>0$ at every $y<1$ and returns ATP levels to $y=1$ at the cost of $u^*\rightarrow\infty$ as $y\rightarrow0$. Each truncation $u_N = c+v_N$ is a bounded polynomial, thus it decreases to below $c/\gamma$ near $y=0$. Since LQR is a bounded affine feedback, while the holding activity $c/\gamma(y)$ grows without bound as $y\rightarrow0$, LQR cannot recover ATP from low levels.

\section{Conclusion}

Al'brekht's method solves one Riccati equation and then an endless sequence
of linear equations. We asked what these equations add up to. For scalar
analytic input-affine systems, Carleman linearization gave the answer. They are the
entries of a single infinite-dimensional Riccati equation in a
lower-triangular Toeplitz matrix. The familiar Riccati equation of the
standard linearization at the equilibrium is its first entry. The nonlinear
terms fill the entries below it. The lifted equation has exactly one solution
with positive diagonal entries. That solution is the derivative of the value
function, and its finite sections converge to it exponentially. What is new
is not the solution. In the scalar case,  it has a closed form, and that closed
form served as our benchmark. What is new is the structure the lift reveals:
a Riccati equation that is decided by its first entry and solved exactly, one
entry at a time. After lifting, the equation of linear optimal control and
the equation of nonlinear optimal control are the same equation. Nothing in
this structure depends on the state dimension. The extension to general
nonlinear systems is the natural next step.

\appendix
\section{Appendix}

\allowdisplaybreaks
\subsection{Proof of Theorem \ref{HJBreformulation.theorem}}\label{app:proofThm1}
By direct computation, we have
\begin{align}
\Psi^\top {\bf P}^\top{\bf A} \Psi
&= \sum_{m=2}^\infty \Big(\sum_{i=1}^{m-1} \sum_{j=1}^{m-i} p(j, i) f_{m-i-j+1}\Big) x^m \nonumber\\
&= \sum_{m=2}^\infty \sum_{k=1}^{m-1} f_{k} \pi_{m-k+1}({\bf P}) x^m \nonumber\\
&= \sum_{m=2}^\infty \Big(\sum_{i+j=m} \sum_{k=1}^{i} \pi_{i-k+2}({\bf P}) f_k \delta_{j1} \Big) x^{m} \nonumber\\
&= \Psi^\top \mathbf{\Pi} \bar{\bf A} \Psi, \label{PA.formula}
\end{align}
where $\delta$ is the standard Kronecker delta symbol.
By a similar argument, we obtain
\begin{align}
\Psi^\top\mathbf{P}^\top \big(\mathbf{B}\Psi+{\bf b}_0\big)
&= \sum_{i=1}^\infty \sum_{j=1}^\infty \Big(\sum_{k=1}^{j+1} p(k,i) g_{j-k+1}\Big) x^{i+j} \nonumber\\
&\quad + g_0 \sum_{i=1}^\infty p(1, i) x^i \nonumber\\
&= \sum_{i=1}^\infty \sum_{j=0}^\infty \Big(\sum_{k'=0}^j p(j-k'+1, i) g_{k'}\Big) x^{i+j} \nonumber\\
&= \sum_{k'=0}^\infty g_{k'}\Big( \sum_{i=1}^\infty \sum_{j'=0}^\infty p(j'+1, i) x^{i+j'}\Big) x^{k'} \nonumber\\
&= \sum_{k'=0}^\infty g_{k'} \sum_{m=2}^\infty \pi_m({\bf P}) x^{m-1+k'} \nonumber\\
&= \sum_{i=1}^\infty \sum_{j=1}^i \pi_{i-j+2}({\bf P}) g_{j-1} x^i
 = \Psi^\top \mathbf{\Pi} {\bf g}. \label{PB.formula}
\end{align}
Applying the expressions \eqref{Q.def} and \eqref{P.def} into \eqref{HJB.linearization}, and then using \eqref{PA.formula} and \eqref{PB.formula} proves \eqref{HJB.linearization.invariant}. The conclusion \eqref{eqn:bilinearHJB.invariant} then follows from combining equations \eqref{eqn:bilinearHJB}, \eqref{P.def} and \eqref{PB.formula}. \hfill $\blacksquare$
\subsection{Proof of Corollary \ref{U.cor}}\label{app:proofCor1}
By 
\eqref{APsi.BPsi.eqn}, 
\eqref{PA.formula} and 
\eqref{PB.formula},
\begin{align*}
\Psi^\top(\mathbf{\Pi}\bar{\mathbf{A}}+\bar{\mathbf{A}}^\top\mathbf{\Pi}^\top)\Psi
  &= 2\Psi^\top{\bf P}^\top{\bf A}\Psi = f(x)U(x),\\
\Psi^\top\mathbf{\Pi}\mathbf{g}
  &= \Psi^\top{\bf P}^\top({\bf B}\Psi+{\bf b}_0) = \tfrac12 g(x)U(x).
\end{align*}
Substituting the above two expressions and \eqref{Q.def} into \eqref{HJB.linearization.invariant} and \eqref{eqn:bilinearHJB.invariant} proves both claims in the corollary. \hfill $\blacksquare$

\subsection{Proof of Theorem \ref{HJBcoefficient.thm}}\label{app:proofThm2}
 By  \eqref{Q.def},  \eqref{PA.formula} and \eqref{PB.formula}, we can reformulate the
 algebraic Riccati-type equation  
 \eqref{HJB.linearization.invariant}
 as 
 \begin{align*} & & \sum_{m=0}^\infty 2\sum_{k=0}^{m}
        f_{m-k+1} \pi_{k+2}({\bf P}) x^{m+2} +\sum_{m=0}^\infty q_{m+2} x^{m+2}\nonumber \\
      & & \quad  -\frac{x^2}{r} \Big (\sum_{n=0}^\infty \sum_{j=0}^n 
\pi_{j+2}({\bf P}) g_{n-j} x^n\Big)^2=0.\end{align*}
Comparing the coefficients of $x^{m+2}, m\ge 0$, we obtain 
\begin{align}\label{HJBcoefficient.thm.pf.eq2}
  &  -\frac{1}{r} \sum_{n=0}^{m}
       \sum_{j=0}^n \sum_{j'=0}^{m-n} g_{n-j} g_{m-n-j'} 
\pi_{j+2}({\bf P}) \pi_{j'+2}({\bf P})\nonumber\\
& \quad + 2\sum_{k=0}^{m}
        f_{m-k+1} \pi_{k+2}({\bf P}) +q_{m+2}  =0, \ \ m\ge 0.
\end{align}
Taking $m=0$ in \eqref{HJBcoefficient.thm.pf.eq2} yields the Riccati equation
of the linearization of the original system \eqref{eqn:GeneralNonlinear} at the origin,
\begin{equation}
\label{HJBcoefficient.eq3+}
q_2 + 2\pi_2({\bf P})\, f_1
    - \frac{g_0^2}{r}\big(\pi_2({\bf P})\big)^2 = 0,
\end{equation}
whose unique positive solution is given by
\eqref{HJBcoefficient.thm.eq1}. Taking $m\ge 1$ in  \eqref{HJBcoefficient.thm.pf.eq2} gives
\begin{align}   \label{HJBcoefficient.thm.pf.eq4}
& 
2\big(f_1-r^{-1} g_0^2 \pi_2({\bf P})\big) \pi_{m+2}({\bf P})\nonumber\\
= & \frac{1}{r} \sum_{n=1}^{m-1}
       \sum_{j=0}^n \sum_{j'=0}^{m-n} g_{n-j} g_{m-n-j'} 
\pi_{j+2}({\bf P}) \pi_{j'+2}({\bf P})\nonumber\\
&  \ \  +\frac{2g_0}{r} \pi_2({\bf P})\sum_{j=0}^{m-1} g_{m-j} \pi_{j+2}({\bf P}) \nonumber\\
& \ \ -2\sum_{k=0}^{m-1} f_{m-k+1} \pi_{k+2}({\bf P})-q_{m+2}.
\end{align}
Then the inductive definition in \eqref{HJBcoefficient.thm.eq2} for  $\pi_{m+2}({\bf P})$, $m\ge 1$, is obtained by dividing by $-2\sqrt{\Delta}$ on both sides of \eqref{HJBcoefficient.thm.pf.eq4} and observing that $f_1 -r^{-1} g_0^2\,\pi_2({\bf P})= -\sqrt{\Delta} < 0$ by  \eqref{HJBcoefficient.thm.eq1}. Since $-2\sqrt{\Delta}\neq0$, the solution is unique. Conversely, for $\pi_m(\mathbf P)=v_{m-1}/2$ in \eqref{Pi.def} with $v_k$ as in \eqref{Vx.expansion}, the series $\sum_{m=2}^\infty\pi_m(\mathbf{P})x^m$ converges by \eqref{vk.def}, has $\pi_2(\mathbf P)>0$ by \eqref{eqn:V(x)}, and satisfies \eqref{HJB.linearization.invariant} by the computation in the proof of Corollary \ref{U.cor}, since $V_x$ solves \eqref{eqn:nonlinearHJB}. Hence, by the uniqueness of this solution, $2\pi_m(\mathbf P)=v_{m-1}$ for all $m\ge2$. \hfill $\blacksquare$

\subsection{Proof of Theorem \ref{finitesection.thm}}\label{app:proofThm3}
Let $v_k, k\ge 1$, be Maclaurin coefficients of the derivative
$V_x$ of the value function $V$  in \eqref{Vx.expansion}.
By Corollary \ref{U.cor}, $U$ in \eqref{U.def} solves \eqref{eqn:nonlinearHJB} and so does $V_x$ in \eqref{eqn:V(x)}. Since $2\pi_2(\mathbf{P})=v_1$ by \eqref{HJBcoefficient.thm.eq1} and \eqref{eqn:V(x)}, and the level equations \eqref{HJBcoefficient.thm.pf.eq2} fix all later coefficients uniquely, we see that
\begin{equation}
\label{finitesection.thm.pfeq1} 
2\pi_m({\bf P})= v_{m-1}, \ m\ge 2.
\end{equation}
This together with \eqref{vk.def}, \eqref{PiN.def} and \eqref{UN.def} implies that
\begin{align*}
|U_N(x)-V_x| & = \Big|\sum_{k=N+1}^\infty v_k x^k\Big|
\le M \sum_{k=N+1}^\infty (|x|/\rho)^k\nonumber\\
& = M (1-|x|/\rho)^{-1} (|x|/\rho)^{N+1}, 
\end{align*}
and the desired estimate \eqref{finitesection.thm.eq1} follows.  \hfill $\blacksquare$







\bibliographystyle{IEEEtran}
\bibliography{carleman}

\end{document}